\documentclass[epsf]{article}
\usepackage{amsfonts}
\usepackage{subfigure} 
\usepackage{epsfig}
\usepackage{amsmath, amsthm}
\newcommand{\be}{\begin{eqnarray}}     	\newcommand{\ee}{\end{eqnarray}}

\newcommand{\vol}{\mathrm{Vol}}

\newcommand{\rem}{\mathrm{Rm}}
\newcommand{\ric}{\mathrm{Ric}}

\title{Curvature tensor under the Ricci flow}
\author{Nata\v{s}a \v{S}e\v{s}um}
\date{}
\theoremstyle{plain}
\newtheorem{dummy}{Dummy}

\theoremstyle{definition}
\newtheorem{corollary}[dummy]{Corollary}

\newtheorem{lemma}[dummy]{Lemma}
\newtheorem{theorem}[dummy]{Theorem}

\newtheorem{definition}[dummy]{Definition}

\newtheorem{claim}[dummy]{Claim}

\begin{document}

\maketitle

\begin{abstract}
Consider the unnormalized Ricci flow $(g_{ij})_t = -2R_{ij}$ for $t\in
[0,T)$, where $T < \infty$. Richard Hamilton showed that if the
curvature operator is uniformly bounded under the flow for all times
$t\in [0,T)$ then the solution can be extended beyond $T$. We prove
that if the Ricci curvature is uniformly bounded under the flow for
all times $t\in [0,T)$, then the curvature tensor has to be
uniformly bounded as well. 
\end{abstract}

\begin{section}{Introduction}

The Ricci flow equation is the evolution equation $(g_{ij})_t =
-2R_{ij}$ introduced by Richard Hamilton in his seminal paper
\cite{hamilton1982}. Short time existence for solutions to the Ricci flow 
on a compact manifold waswas first showen in \cite{hamilton1982},
using the Nash-Moser theorem. Shortly after that De Turck
(\cite{turck}) showed the same thing by modifying the flow by a
reparametrization using a fixed background metric to break the symmetry
that comes from the fact that the Ricci tensor is invariant under the
whole diffeomorphism group of a manifold. 

The Ricci flow equation is a weakly parabolic equation and many nice
regularity theorems have already been proved. A very nice corollary of
the regularity of the Ricci flow is a result on the maximal existence
time for a solution, proved by R. Hamilton in
\cite{hamiltonMA}.

\begin{theorem}[Hamilton]
\label{theorem-theorem_max_time}
For any smooth initial metric on a compact manifold there exists a
maximal time $T$ on which there is a unique smooth solution to the
Ricci flow for $0 \le t < T$. Either $T = \infty$ or the curvature is
unbounded as $t\to T$.
\end{theorem}

The main theorem in this paper is the following 

\begin{theorem}
\label{theorem-theorem_blow_up}
Let $(g_{ij})_t = -2R_{ij}$, for $t\in [0,T)$ be a Ricci flow on a
compact manifold $M$ with $T < \infty$ and with uniformly bounded
Ricci curvatures along the flow. Then the curvature tensor stays
uniformly bounded along the flow.
\end{theorem}

The proof of the theorem uses Perelman's noncollapsing theorem for the
unnormalized Ricci flow that has been proved in \cite{perelman2002}.

Acknowledgements: The author would like to
thank her advisor Gang Tian for suggesting this problem and for all
his help, support and suggestions.

\end{section}

\begin{section}{Background}

$M$ will always denote a compact manifold, and $(g_{ij})_t = -2R_{ij}$
is the unnormalized Ricci flow on $M$. Perelman introduced the
following functional in \cite{perelman2002}

$$W(g,f,\tau) = (4\pi\tau)^{-\frac{n}{2}}\int_M e^{-f}[\tau(|\nabla
f|^2 + R) + f - n] dV_g.$$

Perelman has showed that $W$ is increasing along the flow. A very nice
application of the monotonicity formula for $W$ is Perelman's
noncollapsing theorem for the unnormalized Ricci flow.

\begin{definition}
Let $g_{ij}(t)$ be a smooth solution to the Ricci flow
$(g_{ij})_t=-2R_{ij}(t)$ on $[0,T)$. We say that $g_{ij}(t)$ is
loacally collapsing at $T$, if there is a sequence of times $t_k\to T$
and a sequence of metric balls $B_k = B(p_k,r_k)$ at times $t_k$, such
that $\frac{r_k2}{t_k}$ is bounded, $|\rem|(g_{ij}(t_k))\le r_k^{-2}$ in
$B_k$ and $r_k^{-n}\vol(B_k)\to 0$.
\end{definition}

\begin{theorem}
\label{theorem-perelman_theorem}
If $M$ is closed and $T < \infty$, then $g_{ij}(t)$ is not locally
collapsing at $T$.
\end{theorem}

The corollary of the theorem above states

\begin{corollary}
\label{corollary-corollary_ancient}
Let $g_{ij}(t)$, $t\in [0,T)$ be a solution to the Ricci flow on a
closed manifold $M$, where $T < \infty$. Assume that for some
sequences $t_k\to T$, $p_k\in M$ and some constant $C$ we have $Q_k =
|\rem|(x,t) \le CQ_k$, whenever $t < t_k$. Then a subsequence of
scalings of $g_{ij}(t_k)$ at $p_k$ with factors $Q_k$ converges to a
complete ancient solution to the Ricci flow, which is
$\kappa$-noncollapsed on all scales for some $\kappa > 0$.
\end{corollary}

Before we start proving our theorem, we will mention some estimates
for ball sizes under evolution by the Ricci flow derived by
Glickenstein in \cite{glickenstein2002}, since we will use them in the
proof of our theorem \ref{theorem-theorem_blow_up}. His estimates are
improvements of Hamilton's work on distance comparison at different
times under Ricci flow that can be found in \cite{hamilton1999}.

\begin{lemma}
\label{lemma-lemma_hamilton}
Suppose $(M,g(t)$ are solutions to the Ricci flow for $t\in [0,T)$
such that the Ricci curvatures are uniformly bounded by some constant
$C$. Then for all $\delta > 0$ there exists an $\eta > 0$ such that if
$|t - t_0| < \eta$ then

$$|d_{g(t)}(q,q') - d_{g(t_0)}(q,q')| \le \delta^{\frac{1}{2}}
d_{g(t_0)}(q,q'),$$
for all $q$, $q' \in M$.
\end{lemma}

One can choose $\delta = e^{2C|t - t_0|} - 1$ and $\eta = \ln
\frac{(\delta + 1)}{2C}$.

\begin{lemma}[Glickenstein]
\label{lemma-lemma_size_balls}
If $g(t)$ is a solution to the Ricci flow such that the Ricci
curvature is bounded by constant $1$ then for all $\rho > 0$

$$B_{g(t)}(0, r(t)\rho) \subset B_{g(0)}(0,\rho),$$
$$B_{g(0)}(0, r(t)\rho) \subset B_{g(t)}(0,\rho),$$
where $r(t) = \frac{1}{1 + (e^{2t} -1)^{\frac{1}{2}}}$.
\end{lemma}

The proof of lemma \ref{lemma-lemma_size_balls} uses lemma
\ref{lemma-lemma_hamilton} and triangle inequality for distances.

\end{section}

\begin{section}{Blow up argument}

In this section we will prove theorem (\ref{theorem-theorem_blow_up}).

\begin{proof}

Assume that under the assumtions of our theorem the curvature blows up
at a finite time $T$. That means there exist sequences $t_i\to T$ and
$p_i\in M$, such that 

$$Q_i = |\rem|(p_i,t_i) \ge C^{-1}\max_{M\times [0,t_i]} |\rem|(x,t),$$
and $Q_i \to \infty$ as $i\to \infty$. By corollary
\ref{corollary-corollary_ancient} the scaled metrics $g_i(t) =
Q_ig(t_i + \frac{t}{Q_i})$ converge to a complete ancient solution to
the Ricci flow, $\kappa$-noncollapsed on all scales for some $\kappa >
0$. We have that $\{(M,g_i(t),p_i)\}\to (N,\bar{g}(t),p)$ in a pointed
Gromov-Hausdorff metric, for all $t\in (-\infty, a)$ where $a >
0$. Since the Ricci curvatures of our original flow are uniformly
bounded, we have that $|R(t)|\le C$, where $R(t)$ is a scalar
curvature of $g(t)$. Since $R_i(t) = \frac{R(t)}{Q_i}$ for all $t\in
[-t_iQ_i, (T - t_i)Q_i)$, we get that $R_i(t)\to 0$ as $i\to\infty$,
i.e. our ancient solution $\bar{g}(t)$ has zero scalar curvature for
every $t\in (-\infty,a)$, where $a > 0$. The evolution equation for a
scalar curvature $\bar{R}$ is:

$$\frac{d}{dt}\bar{R} = \Delta\bar{R} + 2|\overline{\ric}|^2.$$

Since $\bar{R}(t)\equiv 0$ for all $t\in(-\infty, a)$ the evolution
equation for $\bar{R}$ implies that $\overline{\ric}(t) = 0$, for all
$t$, i.e. our solution $\bar{g}$ is stationary. Therefore it can be
extended for all times $t\in (-\infty,+\infty)$ to an eternal
solution, where $\bar{g}(t) = \bar{g}$.

Take any $r > 0$. Then:

\begin{equation}
\label{equation-equation_volume_limit}
\frac{\vol B(p,r)}{r^n} = \lim_{i\to\infty}\frac{\vol_i
B_i(p_i,r)}{r^n},
\end{equation}
where the volume and the ball $B(p,r)$ on the LHS of
\ref{equation-equation_volume_limit} are considered in metric
$\bar{g}$, while on the RHS of \ref{equation-equation_volume_limit} we
consider metric $g_i(0) = Q_ig(t_i)$. Furthermore,

$$\frac{\vol B(p,r)}{r^n} =
\lim_{i\to\infty}\frac{\vol_{g(t_i)}B_{g(t_i)}(p_i,rQ_i^{-\frac{1}{2}})}
{(rQ_i^{-\frac{1}{2}})^n}.$$

Since the Ricci curvatures of our original metrics $g(t)$ are
uniformly bounded, from the evolution equation for $\vol_t$ we get
that $\vol_{g(t)}M \le \tilde{C}$ for all $t\in [0,T)$, for some
constant $\tilde{C}$ that does not depend on $t$. Take any $\epsilon >
0$. Choose $i_0$ such that for all $i \ge i_0$ we have the following
estimates:

\begin{equation}
\label{equation-equation_choice_i1}
e^{-\frac{C}{\sqrt{n}}}|t_i - t_{i_0}| > (1 - \frac{\epsilon}{2}),
\end{equation}

\begin{equation}
\label{equation-equation_choice_i2}
(\frac{1}{1 + (e^{2C(t_i-t_{i_0})} - 1)^{\frac{1}{2}}})^n > 1 -
\frac{\epsilon}{2},
\end{equation}

and

\begin{equation}
\label{equation-equation_epsilon}
\frac{\vol B(p,r)}{r^n} \ge
\frac{\vol_{g(t_i)}B_{g(t_i)}(p_i,rQ_i^{\frac{-1}{2}})}
{(rQ_i^{-\frac{1}{2}})^n} - \frac{\epsilon}{2},
\end{equation}
hold for all $i \ge i_0$, where $C$ is a constant in the statement of
our theorem \ref{theorem-theorem_blow_up}. We can choose such $i_0$
because of \ref{equation-equation_volume_limit} and the fact that
$t_i\to T$ as $i\to \infty$.

Let $r_i = \frac{1}{1 + (e^{2C(t_i - t_{i_0})} - 1)^{\frac{1}{2}}}$.
Lemma \ref{lemma-lemma_size_balls} tells us how the size of balls
change with the Ricci flow. Since $t_i \ge t_{i_0}$, we have that
$B_{g(t_{i_0})}(p_i,rr_iQ_i^{-\frac{1}{2}}) \subset
B_{g(t_i)}(p_i,rQ_i^{-\frac{1}{2}})$. From this observation and
estimate \ref{equation-equation_epsilon} we have:

$$\frac{\vol B(p,r)}{r^n} \ge \frac{\vol_{g(t_i)}B_{g(t_{i_0})}(p_i,
rr_iQ_i^{-\frac{1}{2}})}{(r\sqrt{Q_i})^n} - \frac{\epsilon}{2}.$$

Denote by $\tilde{r}_i = rr_iQ_i^{-\frac{1}{2}}$.  From the evolution
equation for $\vol_t$ we get

\begin{equation}
\label{equation-equation_main_contra1}
\vol_{t_i} = \vol_{t_{i_0}}e^{-\int_{t_{i_0}}^{t_i} R du}.
\end{equation}
Integrating equation (\ref{equation-equation_main_contra1}) over a ball
$B_{g(t_{i_0}}(p_i,\tilde{r}_i)$ gives

\begin{eqnarray}
\label{eqiation-equation_main_contra}
\vol_{g(t_i)}B_{g(t_{i_0})}(p_i,\tilde{r}_i) &=& 
\int_{B_{g(t_{i_0})}} e^{-\int_{t_{i_0}}^{t_i} R du} dV_{g(t_{i_0})} 
\nonumber \\
&\ge& \vol_{g(t_{i_0})}B_{g(t_{i_0})}(p_i,\tilde{r}_i)
e^{-|t_i - t_{i_0}|\frac{C}{\sqrt{n}}}
\nonumber \\
&\ge& (1 - \frac{\epsilon}{2})\vol_{g(t_{i_0})}
B_{g(t_{i_0})}(p_i,\tilde{r}_i),
\end{eqnarray}
where we have used estimate \ref{equation-equation_choice_i1} and the
fact that $|R|^2 \le \frac{|\ric|^2}{n} \le \frac{C^2}{n}$. This gives
the following estimate

\begin{equation}
\label{equation-equation_epsilon_estimate}
\frac{\vol B(p,r)}{r^n} \ge (1 -
\frac{\epsilon}{2})\frac{\vol_{g(t_{i_0})}B_{g(t_{i_0})}(p_i,
\tilde{r}_i)}{(r\sqrt{Q_i})^n} - \frac{\epsilon}{2}.
\end{equation}

We have a fixed metric $g(t_{i_0})$ on the RHS of the inequality
\ref{equation-equation_epsilon_estimate}. Our manifold $M$ is compact
and therefore we have a uniform asymptotic volume expansion:

\begin{equation}
\label{equation-equation_asymptotic}
\vol_{t_{i_0}}B_{t_{i_0}}(p_i, \tilde{r}_i) = \tilde{r}_i^n\omega_n (1
- \frac{R(g(t_{i_0})) r^2}{6(n+2)Q_i(1 + (e^{2C(t_i - t_{i_0})} -
1)^{\frac{1}{2}})^2} + o(\tilde{r}_i^2)),
\end{equation}
where $\omega_n$ is a volume of a unit euclidean ball.

Combining the estimate \ref{equation-equation_epsilon_estimate} with
the asymptotic expansion \ref{equation-equation_asymptotic} and
letting $i\to\infty$ we get:

$$\frac{\vol B(p,r)}{r^n} \ge r_i^n\omega_n (1 - \frac{\epsilon}{2}) -
\epsilon.$$
By estimate \ref{equation-equation_choice_i2} we have that $r_i^n > 1
- \frac{\epsilon}{2}$ for all $i \ge i_0$ and therefore

$$\frac{\vol B(p,r)}{r^n} \ge \omega^n (1 - \frac{\epsilon}{2})^2 -
\epsilon.$$ 
Since $\epsilon > 0$ was arbitrary, we have that
$\frac{\vol B(p,r)}{r^n} \ge \omega_n$.

On the other hand Bishop-Gromov volume comparison principle applied to
a Ricci flat, complete manifold $(N,\bar{g})$ gives 

$$\frac{\vol B(p,r)}{r^n} \le \frac{\vol B(p, \delta)}{\delta^n},$$
for all $\delta \le r$. When $\delta\to 0$, the RHS of the previous
inequality tend to $\omega_n$. Therefore for every $r > 0$ we would
have that $\frac{\vol B(p,r)}{r^n} = w_n$. Ricci flat, complete
manifold with this property has to be the Euclidean space and therefore
$\overline{\rem} \equiv 0$ on $N$. This contradicts the fact that
$\overline{\rem}(p,0) = 1$.

Therefore, a curvature of an unormalized flow can not blow up in a
finite time when the Ricci curvatures are uniformly bounded for all
times for which a solution of a Ricci flow exists.

\end{proof}

We have immediatelly the following simple corollary of theorem
\ref{theorem-theorem_blow_up}.

\begin{corollary}
Let $g(t)$ be a solution to $(g_{ij})_t = -2R_{ij}$ with $|\ric| \le
C$ uniformly for all times when the solution exist. Then the solution
exists for all times $t\in [0,\infty)$.
\end{corollary}

\begin{proof}

This is a simple consequence of theorems \ref{theorem-theorem_blow_up}
and \ref{theorem-theorem_max_time}.

\end{proof}

\begin{corollary}
Let $(g_{ij})_t = -2R_{ij} + \frac{1}{2\tau}g_{ij}$ be a flow on a
closed manifold $M$, which solution exists for $t\in [0,T)$, where $T
< \infty$. If $|\ric| \le C$ for all $t\in [0,T)$, then the solution
can be extended past time $T$.
\end{corollary}

\begin{proof}
Let $\bar{g}(s) = g(t(s))c(s)$, where $c(s) = 1 - \frac{s}{\tau}$ and
$t(s) = -\tau\ln(1 - \frac{s}{\tau})$. Then $|\ric|_{\bar{g}(s)} \le
\frac{C}{c(s)} \le Ce^{\frac{T}{\tau}}$ for all $s\in [0, \tau(1 -
e^{-\frac{T}{\tau}}))$. We can apply theorem
(\ref{theorem-theorem_blow_up}) to the unnormalized flow $\bar{g}(s)$
to get that it can be extended past time $\tau(1 -
e^{-\frac{T}{\tau}}))$. If we go back to flow $g(t)$, it means that
$g(t)$ can be extended past time $T$.
\end{proof}

Since Perelman's noncollapsing theorem plays an important role in a
study of the Ricci flow, we would like to say few words about it in
the case of normalized flow.

\begin{claim}
\label{claim-claim_translation}
Let $(g_{ij})_t = -2R_{ij} + \frac{1}{\tau}g_{ij}$ for $t\in
[0,\infty)$ on a closed manifold $M$. Then $g_{ij}(t)$ is not locally
collapsing at $\infty$.
\end{claim}

Before we start with the proof of the claim, we will write down the
definition for collapsing at the $\infty$.

\begin{definition}
Let $g_{ij}(t)$ be a solution to $(g_{ij})_t = -2R_{ij} +
\frac{1}{\tau}g_{ij}$. We say that $g_{ij}(t)$ is locally collapsing
at $\infty$, if there is a sequence of times $t_i \to \infty$ and a
sequence of metric balls $B_k = B(p_k,r_k)$ at times $t_k$, such that
$\frac{r_k^2}{t_k}$ is bounded, $|\rem|g_{ij}(t_k))| \le r_k^{-2}$ and
$r_k^{-n}\vol(B_k) \to 0$.
\end{definition}

Notice that in the proof of \ref{theorem-perelman_theorem} we need
just that $|\ric|(g_{ij}(t_k))| \le r_k^{-2}$.

\begin{proof}[proof of claim \ref{claim-claim_translation}]

Assume that there exists a sequence of collapsing balls $B_k =
B(p_k,r_k)$ at times $t_k\to\infty$ where $\frac{r_k^2}{t_k}$ is
bounded. Let $\tilde{g}(s) = c(s)g(t(s))$, where $c(s) = 1 -
\frac{s}{\tau}$ and $t(s) = -\tau\ln(1 - \frac{s}{\tau})$. This choice
of $c(s)$ and $t(s)$ gives us that $\tilde{g}(t)$ satisfies
$(\tilde{g}_{ij})_t = -2R_{ij}$, for $s\in [0,\tau)$. Perelman's
noncollapsing theorem \ref{theorem-perelman_theorem} gives us that
$\tilde{g}(s)$ is not locally collapsing at $\tau$. Let $s_k$ be a
sequence such that $t(s_k) = t_k$. By our assumption there exists a
sequence $\epsilon_k \to 0$ as $k\to\infty$ such that

$$\vol_{g}B_g(p_k,r_k) \le \epsilon_k r_k^n,$$
at time $t_k$. This we can write as

$$\vol_{\tilde{g}}B_{\tilde{g}}(p_k,\sqrt{c(s_k)}r_k) \le
\epsilon_k(r_k\sqrt{c(s_k)})^n,$$
at time $s_k$. $\frac{c(s_k)r_k^2}{s_k^2}$ is bounded ( for
$\frac{c(s_k)r_k^2}{s_k^2} = (-\frac{\tau\ln(1 -
\frac{s_k}{\tau})\cdot (1 - \frac{s_k}{\tau})}{s_k^2})
(\frac{r_k^2}{t_k})$, where the first bracket is bounded since $s_k\to
1$ and $\ln(1 - s_k) (1 - s_k)\to 0$ as $k\to\infty$ and the second
bracket is bounded by assumption).

This means a sequence of balls $B(p_k,\sqrt{c(s_k)}r_k)$ at times
$s_k$ is a sequence of collapsing balls at $\tau$ for flow
$\tilde{g}_{ij}(s)$, which contradicts theorem
\ref{theorem-perelman_theorem}.

\end{proof}

\begin{lemma}
\label{lemma-lemma_consequence}

Let $g_{ij}(t)$, $t\in [0,\infty)$ be a solution to the Ricci flow
$(g_{ij})_t = -2R_{ij} + \frac{1}{\tau}g_{ij}$ on a closed manifold
$M$, with uniformly bounded Ricci curvatures. Assume that for some
sequences $t_k\to\infty$, $p_k\in M$ and some constant $C$ we have
$Q_k = |\rem|(p_k,t_k) \to\infty$ and $|\rem|(x,t) \le CQ_k$ when $t
\le t_k$. Then a subsequence of scalings of $g_{ij}(t_k)$ at $p_k$
with factors $Q_k$ converges to a complete, eternal solution to the
Ricci flow that is noncollapsed on all scales for some $\kappa > 0$.

\end{lemma}

\begin{proof}

This lemma is a simple corollary of claim
\ref{claim-claim_translation}. The proof is the same as a proof of a
corollary of Perelman's noncollapsing theorem in section $4$ of paper
\cite{perelman2002}.

\end{proof}

\end{section}

\begin{section}{The volume form under the Ricci flow in dimension three}

In this section we would like to mention one nice application of the
proof of theorem \ref{theorem-theorem_blow_up} to $3$ dimensional
manifolds, observed by Richard Hamilton.

\begin{theorem}
\label{theorem-therem_volume_element}
Fix $t_0\in [0,T)$.  Let $(g_{ij})_t = -2R_{ij}$ be the Ricci flow on
a 3-dimensional compact manifold, for $t\in [0,T)$, where $T <
\infty$. If $g(t)$ is singular at time $T$, then $\lim_{t\to
T}\hat{\vol}_t = 0$, where $\hat{\vol}_t = \min_{x\in M}
\frac{\vol_t}{\vol_{t_0}}$.
\end{theorem}

\begin{proof}

We will prove theorem \ref{theorem-therem_volume_element} by
contradiction. Assume that there exist a sequence of times $t_i \to T$
and $\delta > 0$ so that $\hat{\vol}_{t_i} > \delta$ for all $i$. That
would imply $\frac{\vol_{t_i}}{\vol_{t_0}} > \delta$ for all $i$ and
all $x\in M$.

\begin{claim}
There exists $C = C(g(0), T)$ so that $\frac{\vol_t}{\vol_{t_0}} \le
C$ for all $t\in [0,T)$ and all $x\in M$.
\end{claim}

\begin{proof}

The evolution equation for $\ln \vol_t$ is

\begin{equation}
\label{equation-equation_volume}
\frac{d}{dt}\ln\vol_t = -R.
\end{equation}

The evolution equation for the scalar cuvature $R(t)$

$$\frac{d}{dt}R = \Delta R + 2|\ric|^2,$$ 
implies by a starightforward maximum principle argument that 
at any time $t\in [0,T)$
$$R(t) \ge \frac{1}{(\min R(0))^{-1} - 2t/3} \ge -\frac{3}{2(t +
C)},$$ for some constant $C$. 
If we integrate the equation
(\ref{equation-equation_volume}) over $s\in [t_0,t]$ for any 
$t\in [t_0,T)$ we will get

\begin{eqnarray*}
\ln\vol_t - \ln\vol_{t_0} &=& -\int_{t_0}^t R(s) ds \\
&\le& \int_0^T \frac{3}{2(s + C)} ds \le C_1.
\end{eqnarray*}
Now we easily get the statement of the claim.

\end{proof}

From the equation (\ref{equation-equation_volume}) we get

$$\ln\vol_{t_i} - \ln\vol_{t_0} = -\int_{t_0}^{t_i}R(s) ds.$$ Since
$\ln\delta < \ln\frac{\vol_{t_i}}{\vol_{t_0}} \le \ln C$, where $C$ is
a constant from the claim above, we get

$$|\int_{t_0}^T R(s) ds| \le C,$$ 
where $C$ may be some other constant
depending on $g(t_0)$ and $T$, but nonetheless we will denote all such
constants by the same symbol.

In dimension $3$ we have a pinching estimate, i.e. there exists
constants $C_1$ and $C_2$ so that

$$|\ric| \le C_1(R + C_2),$$ for all times $t\in [0,T)$. This gives
that $\int_{t_0}^T |\ric| ds \le C_1\int_{t_0}^T (R(s) + C_2) ds \le
C$ and the esrimate does not depend on $x\in M$.

If $g(t)$ has a singularity at time $T$, by theorem
\ref{theorem-theorem_blow_up} we know that there exist a sequence of
times $s_i\to T$, sequence of points $x_i\in M$ and a constant $C$ so
that $Q_i = |\ric|(x_i,s_i) \ge
C^{-1}\max_{M\times[0,s_i]}|\ric|(x,t)$ and $Q_i\to\infty$. In
dimension $3$ the curvature tensor is controled by the Ricci
curvature. Let $g_i(s) = Q_ig(s_i + sQ_i^{-1})$. By Perelman's
noncollapsing theorem and Hamilton's compactness theorem there exists
a subsequence $g_i$ so that $(M, g_i(s), x_i)\to (N,\bar{g}(s), x)$,
where $(N,\bar{g}(s))$ is a complete, non-flat, ancient solution to
the Ricci flow. When the Ricci curvature blows up, that happens with
the scalar curvature as well because of the pinching estimate in
dimension $3$. Because of the pinching estimate, when the scalar
curvature $R$ is big, any negative curvature is small compared to the
positive ones. That is why after rescaling by a factor $Q_i$ and after
taking a limit when $i\to\infty$ we get that our ancient, complete,
non-flat solution $(N, \bar{g}(s))$ has bounded, non-negative
sectional curvatures. This allows us to apply the relative volume
comparison theorem to $(N,\bar{g}(s))$ to conclude that for every $r >
0$

$$\frac{\vol_{\bar{g}(0)} B(x,r)}{r^n} \le w_n,$$
where $w_n$ is the volume of a unit euclidean ball. 
Fix $\epsilon > 0$. We can show that 

$$\frac{\vol_{\bar{g}(0)} B(x,r)}{r^n} \ge w_n - \epsilon,$$ in the
same way as in theorem \ref{theorem-theorem_blow_up}, since we used
only integral bound on $|\ric|$ in order to control changes in volumes
and distances for $t$ close to $T$. Since this estimate holds for
every $\epsilon > 0$, we can conclude that $(N, \bar{g}(0))$ would
have to be the euclidean space, which is not possible since
$|\ric|(x,0) = 1$.

\end{proof}

\end{section}

\end{document}